\documentclass{amsart}
\usepackage{amssymb} 
\usepackage{amsthm}
\usepackage{graphicx}
\usepackage{amscd}
\usepackage{color}
\usepackage{setspace}
\usepackage[all]{xy}
 \usepackage[colorlinks=true]{hyperref}
 \hypersetup{urlcolor=blue, citecolor=red}
\theoremstyle{plain}
\newtheorem{theorem}{Theorem}
\newtheorem{lemma}{Lemma}
\newtheorem{corollary}{Corollary}

\theoremstyle{definition}
\newtheorem{definition}{Definition}

\newtheorem{example}{Example}

\theoremstyle{remark}

\newtheorem*{remark}{Remark}

\newcommand{\N}{\mathbb{N}}
\newcommand{\Q}{\mathbb{Q}}
\newcommand{\R}{\mathbb{R}}
\newcommand{\Z}{\mathbb{Z}}

\newcommand{\A}{\alpha}

\newcommand{\abs}[1]{\left|#1\right|}

\newcommand{\integer}[1]{\left[ #1 \right]}

\newcommand{\ST}{\textrm{ }\mathbf{:}\textrm{ }}

\title[Heaviness in Symbolic Dynamics]{Heaviness in Symbolic Dynamics:\\Substitution and Sturmian Systems}
\author[David Ralston]{}
\email{ralston@math.ohio-state.edu}
\begin{document}
\maketitle

\centerline{\scshape David Ralston }
\medskip
{\footnotesize
\centerline{Department of Mathematics, The Ohio State University}
\centerline{231 W. 18th Avenue, Columbus, OH  43210, USA}
} 

\bigskip

\renewcommand{\labelenumi}{\Roman{enumi}.}
\renewcommand{\theenumi}{\Roman{enumi}}

\begin{abstract}
\textit{Heaviness} refers to a sequence of partial sums maintaining a certain lower bound and was recently introduced and studied in \cite{ralston1}.  After a review of basic properties to familiarize the reader with the ideas of heaviness, general principles of heaviness in symbolic dynamics are introduced.  The classical Morse sequence is used to study a specific example of heaviness in a system with nontrivial rational eigenvalues.  To contrast, Sturmian sequences are examined, including a new condition for a sequence to be Sturmian.
\end{abstract}

\maketitle

\section{Introduction}

Dynamical systems devotes much attention to the asymptotic behavior of points or other elements of a system.  While asymptotic properties are extremely important, they are in a sense not observable; an observer monitoring a closed system can only ever observe a finite window of time.  Suppose that an observer is capable of monitoring the output of a function $f$ over a finite portion of an orbit $x, T(x), T^2(x), \ldots, T^n(x)$, and keeps a record of the associated partial sums.  Finite observations do not lend themselves to discussion of limits, but any observer might be concerned with \emph{extremal} behavior of the partial sums (a motivation similar to, but distinct from, that in the study of large deviations).  With this restriction and motivation in mind, we define the \emph{heavy set} (subject to various restrictions to be outlined later) to be those points in a system for which these partial sums maintain a natural lower bound over a natural collection of finite ranges.

In applying these notions to symbolic dynamics, we will note a distinction between systems with rational eigenvalues (\S\ref{section-substitution}) and a class of systems with no rational eigenvalues (\S\ref{section-sturmian}).  Specifically, the Morse sequence defines a nontrivial system with an abundance of heavy points, while Sturmian sequences have a scarcity of heavy points.  Furthermore, Sturmian sequences are most frequently defined with a restriction on the allowable weights of subwords, and heaviness will be concerned with establishing bounds on weights of words, so a connection between the two ideas is developed, most significantly in Theorem \ref{theorem-sturmian equivalence}.

\subsection{Background and terminology}
Before defining heaviness formally in \S\ref{section-heaviness defined section}, it is necessary to establish our framework and notation.  Let $\{\Omega, \mu\}$ be a probability measure space ($\mu(\Omega)=1$).  If $T:\Omega \rightarrow \Omega$ is $\mu$-measurable, and $\mu(T^{-1} \Gamma)=\mu(\Gamma)$ for all $\mu$-measurable $\Gamma \subset \Omega$, then we say that $T$ \emph{preserves $\mu$}, and $\{\Omega, \mu, T\}$ is a \emph{probability measure preserving system}.  In this situation, let $f\in L^1(\Omega,\mu)$.  If the only functions $f$ such that $f\circ T = f$ almost-everywhere are themselves almost-everywhere constant, then $T$ is called \emph{ergodic}, and if $\mu$ is the only preserved probability measure, $T$ is called \emph{uniquely ergodic}.  For $\omega \in \Omega$, define $S_n(\omega)$ recursively: $S_0(\omega)=0$, $S_{n+1}(\omega)=S_n(\omega)+f\circ T^n (\omega)$ (if $T$ is invertible, we may use this relation to define for all $n \in \Z$):
\begin{equation}\label{eqn-partial sums}S_n(\omega)=\sum_{i=0}^{n-1}f\circ T^i (\omega)\hspace{.1 in}(n \geq 0),\hspace{.1 in} S_n(\omega)=-\sum_{i}^{-1}f\circ T^{i}(\omega)\hspace{.1 in}(n<0).\end{equation}

In the vein of classical concerns of asymptotic behavior over infinite time frames, define:
\begin{equation}\label{lower star functions} f_*(\omega)=\liminf_{n \rightarrow \infty} \frac{1}{n}S_n(\omega), \end{equation}
noting that the Birkhoff Ergodic Theorem guarantees that $\lim_{n \rightarrow \infty}n^{-1}S_n(\omega)$ exists almost everywhere.

We may use the fact that $T$ preserves $\mu$ to derive $\int_{\Omega}S_n(x)d \mu = n \int_{\Omega} f d\mu$.  So, in line with our model of an observer of finite time periods, we define:
\begin{definition}{The \emph{heavy set for $f$ (relative to $T$) between times $m$ and $n$ ($m,n \in \Z$, $m \leq n$)} is given by: \begin{equation}\label{heavy m,n definition}\mathcal{H}_T^f(m,n)=\left\{\omega \ST \forall i, \hspace{.1 in} m \leq i \leq n, \hspace{.1 in} S_i(\omega) \geq i \int_\Omega f d\mu\right\}.\end{equation}
\noindent In the common event that $m=0$, the set is called \emph{heavy through time $n$}.  We use the shorthand: \[\mathcal{H}_T^f(\N)=\bigcap_{i=0}^{\infty}\mathcal{H}_T^f(0,i),\hspace{.1 in}\mathcal{H}_T^f(\Z)=\bigcap_{i=0}^{\infty}\mathcal{H}_T^f(-i,i),\]
\noindent to define the \emph{heavy sets over $\N$ and $\Z$}.  Any use of negative times requires $T^{-1}$ to exist.}
 \label{heavinessdefinition}
\end{definition}

If $T$ and $f$ are clear from the context, we will simply write $\mathcal{H}(m,n)$, $\mathcal{H}(\N)$, or $\mathcal{H}(\Z)$. These sets represent points whose partial sums meet or exceed the average value of the partial sums over the range of prescribed times.  Given the emphasis on finite time periods in defining heavy sets, it is worth pointing out that stating $x \in \mathcal{H}(\N)$ or $x \in \mathcal{H}(\Z)$ should not be read as a statement about the behavior of $S_n(x)$ on an infinite time frame, but rather about \emph{all finite times}.  We now present a pair of theorems regarding the existence of such points.

\begin{theorem}{If $\{\Omega, \mu, T\}$ is a measure-preserving system, and $f \in L^1(\Omega,\mu)$, then $\mu \left( \mathcal{H}_T^f(0,n) \right)>0$ for any $n \in \N$.  Furthermore, if $T$ is invertible, then $T$ is ergodic if and only if for any $m<n$, $m,n \in \Z$ and $f \in L^1(\Omega,\mu)$, $\mu \left( \mathcal{H}_T^f(m,n) \right)>0$.}\label{positivemeasureheavytheorem}
\begin{proof}
 For a proof of this theorem, see \cite{ralston1}.
\end{proof}
\end{theorem}

\begin{theorem}{If $\{\Omega, \mu, T\}$ is a continuous measure-preserving system on a compact probability space, then $\mathcal{H}_T^f(\N) \neq \emptyset$ for any upper semi-continuous $f \in L^1(\Omega,\mu)$.  Furthermore, if $f$ is continuous and $T$ is invertible and transitive, then $\mathcal{H}_T^f(\Z) \neq \emptyset$.}\label{pereslemmaandextension}
\begin{proof}
 For a proof of this theorem, see \cite{ralston1}.  The first statement follows as a corollary from Theorem \ref{positivemeasureheavytheorem}, but an alternate proof may be found as a lemma of Y. Peres \cite{peres}.
\end{proof}
\end{theorem}
\begin{remark}
 If $\Gamma$ is closed, then $\mathcal{H}_T^{\chi_{\Gamma}}(\N)\neq \emptyset$ (as $\chi_{\Gamma}$ is upper semi-continuous), but $\mathcal{H}_T^{\chi_{\Gamma}}(\Z)\neq \emptyset$ is true in general only if $\Gamma$ is clopen (in this case, $\chi_A$ is continuous).
\end{remark}

\begin{remark}
 It is a common mistake to assume that $\mu \left( \mathcal{H}_T^f(\N) \right) =0$.  This claim is obviously false for functions $f$ which are constant almost everywhere, for any function $f$ on an atomic system (any nonempty set is of positive measure), and in general for nonergodic $T$ (see Corollary \ref{nonergodic corollary}).

 In \S\ref{classicalmorse}, we give an example of a uniquely ergodic system without atoms, and a function $f$ which is not constant, such that $\mu(\mathcal{H}(\N))>0$.  However, as proved in \cite{halasz}, if $T$ is ergodic, then for almost every $\omega \in \Omega$, there exists some $N=N(\omega)$, $0<N<\infty$, for which $S_N(\omega) \leq N \int_{\Omega}f d\mu$.  There is no contradiction between this fact and Theorem \ref{positivemeasureheavytheorem}; if $\mathcal{H}(\N)$ is not a null set, then for almost all $\omega \in \mathcal{H}(\N)$ there is some finite $N>0$ such that $S_N(\omega)=N\int_{\Omega}f d\mu$.  This situation will be investigated in \S\ref{substitution systems}.
\end{remark}

Proceeding to definitions specific to symbolic dynamics, the reader familiar with standard terminology may skip to $\S\ref{section-heaviness defined section}$, except to note the nonstandard Definitions \ref{subword notation} and \ref{reversing operator}.  We let an \emph{alphabet} $\mathcal{A}$ be a subset of $\R$, and elements of $\mathcal{A}$ are called \emph{letters}.  An element $A \in \mathcal{A}^n$ is called a \emph{word} of \emph{length $n$} over $\mathcal{A}$.  We sometimes write $\abs{A}$.  An element $X \in \mathcal{A}^{\N}$ is said to be a \emph{sequence}, and $X \in \mathcal{A}^{\Z}$ is a \emph{bi-sequence}.  In the frequent event that $\mathcal{A}=\{0,1\}$, the word is called \emph{binary}.  The somewhat nonstandard definition of an alphabet to be an arbitrary subset of $\R$ is motivated by the canonical relation of sequences in $\mathcal{A}^{\N}$ to sequences $\{f(T^i \omega)\}_{i=0,1,2,\ldots}$ in a probability measure preserving system $\{\Omega, \mu, T\}$ along with $f \in L^1(\Omega,\mu)$ and $\omega \in \Omega$.  In this case, $\mathcal{A}$ is the range of $f$.  By assuming $\Omega$ to be compact, continuity of $f$ implies compactness of $\mathcal{A}$, and in the common scenario that $f$ is the characteristic function of a set, $\mathcal{A}=\{0,1\}$.  By relating points in $\mathcal{A}^{\N}$ and sequences $\{f(T^i \omega)\}$ ($i=0,1,\ldots$), then, the shift operator on the space $\mathcal{A}^{\N}$  is analogous to the transformation $T:\Omega \rightarrow \Omega$, and heaviness statements about measure preserving systems in general may be interpreted as statements about heaviness in shift systems.

For a binary word, we define the \emph{conjugate} of $A=a_0 a_1 \ldots$, denoted $\overline{A}=\overline{a_0} \hspace{.05 in} \overline{a_1} \ldots$, by setting $\overline{0}=1$ and $\overline{1}=0$.  For any word $A=a_0 \ldots a_{n-1}$ of length $n<\infty$ over any alphabet, the \emph{transpose} of $A$ is denoted and defined by $A^T=a_{n-1} a_{n-2} \ldots a_0$ (so that $A^T_i=A_{n-1-i}$ for $i=0,\ldots,n-1$).  Given two words $A=a_0 \ldots a_{m-1}$, $B=b_0 \ldots b_{n-1}$ of finite lengths $m$ and $n$, the \emph{concatenation} of $A$ and $B$ is the word of length $m+n$ given by $AB=a_0\ldots a_{m-1} b_0 \ldots b_{n-1}$.  Given a word $A \in \mathcal{A}^n$, the \emph{weight} and \emph{average weight}, respectively, are: \[w(A)=\sum_{i=0}^{n-1} a_i, \hspace{.1 in} \overline{w}(A)=\frac{1}{n}w(A).\]

A word $A$ of length $n$ is said to be a \emph{factor} of another word (or sequence) $B$ of length $m\geq n$ if there is some $j \in \N$ so that $a_i=b_{i+j}$ for $i=0,1,\ldots,n-1$.  The \emph{complexity function} for a sequence or bi-sequence $X$ (over a finite alphabet $\mathcal{A}$) is given by \[p(n)=\#\{A \in \mathcal{A}^n \ST A \textrm{ is a factor of $X$}\}.\]  A binary sequence $X$ is said to be of \emph{minimal complexity} if $p(n)+1$ (any sequence $X$ for which $p(n) \leq n$ for some $n$ is eventually periodic - see \cite{coven-hedlund}).  If $A$ is a word of length $n<\infty$ such that there are two distinct letters $\alpha,\beta \in \mathcal{A}$ such that $A\alpha$ and $A\beta$ are both factors of $X$, then $A$ is said to be a \emph{right special factor}.  If there are two distinct letters $\alpha$ and $\beta$ such that $\alpha A$ and $\beta A$ are both factors of $X$, then $A$ is said to be a \emph{left special factor}.

Given a sequence $X \in \mathcal{A}^{\N}$ or $\mathcal{A}^{\Z}$, we define the sequence $\sigma(X)$ by $\sigma (X )_n=X_{n+1}$.  If $\mathcal{A}$ is compact, then so are $\mathcal{A}^{\N}$ and $\mathcal{A}^{\Z}$ (in the product topology), and therefore $\Omega= \overline{\left\{\sigma^n X\right\}_{n=0}^{\infty}}$ is compact. It is seen that $\sigma$ is now a continuous map of a compact space. The \emph{system generated by the sequence $X$} is the topological dynamical system $\left\{ \overline{\{\sigma^n X\}_{n=0}^{\infty}},\sigma\right\}$.

\begin{definition}{\label{subword notation} Let $A=a_0 \ldots a_{n-1}$ and $0\leq i\leq j \leq n$.  Then $A_{i,j}=a_i \ldots a_{j-1}$ is a word of length $j-i$, beginning at index $i$ (note that $A=A_{0,n}$, and $A_{i,i}$ is the empty word of length zero).}
\end{definition}

\begin{definition}{\label{reversing operator} Let $A=a_0 \ldots a_{n-1}$ be a word of length $n$, over alphabet $\mathcal{A}$.  Then the \emph{reversal} of $A$ is the word over the alphabet $-\mathcal{A}=\{-\alpha \ST \alpha \in\mathcal{A}\}$, defined and notated by \[\rho(A)_i=-a_{n-1-i}, \hspace{.1 in}i=0,1,\ldots,n-1.\]  That is, $\rho(A)$ is the transpose of $A$, with a negative sign on all entries.  For $n<0$, $\rho(A_{n,0})=(-a_{-1})(-a_{-2})\ldots(-a_{n})$.  If we define $-A=(-a_0)(-a_1)\ldots(-a_{n-1})$, then $\rho(A)=-(A^T)=(-A)^T$.  For any word $A$, $w(A)=-w(\rho(A))$.  }
\end{definition}
\begin{remark}
 {Compare with \eqref{eqn-partial sums}, our partial sums over negative times.  By defining \mbox{$f:\mathcal{A}^{\Z} \rightarrow \R$} by $f(X)=x_0$, and $\sigma$ is the shift operator, for $n \geq 0$, $S_n(X)=w(X_{0,n})$, and for $n \leq 0$, $S_n(X)=w(\rho(X_{n,0}))$ (for $n=0$, both equal $0$, the weight of the empty word).  This relation is the motivation for defining $\rho$.}
\end{remark}

\section{Heaviness in symbolic dynamics}\label{section-heaviness defined section}
There are two fruitful ways to define heaviness in symbolic dynamics.  The first, \emph{$\A$-heaviness} (\S\ref{alpha heaviness section}), is a direct analogue of the definition of heaviness in Theorem \ref{positivemeasureheavytheorem}.  The second way to view heaviness, \emph{local heaviness} (\S\ref{local heaviness section}), is more combinatorial in nature.  In presenting both views, we will spend some time to familiarize the reader with the definitions by presenting several theorems regarding the existence of such phenomena in very general settings, before we proceed to considering any specific systems.

\subsection{$\A$-Heaviness}\label{alpha heaviness section}
In situations where we are interested in a global target for heaviness, some fixed $\A$ which will act as a lower bound on our partial averages, we proceed as follows:
\begin{definition}{Let $A=a_0 \ldots a_{n-1}$.  $A$ is said to be \emph{$\A$-heavy} (\emph{$\A$-light}) if for all \mbox{$1 \leq i \leq n$}: \begin{equation}\label{alpha heavy word eqn}\overline{w}(A_{0,i})\geq \A \hspace{.2 in} (\overline{w}(A_{0,i})\leq \A).\end{equation}  If $X$ is a sequence, then $X$ is $\A$-heavy ($\A$-light) if for all $i \in \N$: \[\overline{w}(X_{0,i})\geq \A \hspace{.2 in} (\overline{w}(X_{0,i}) \leq \A).\]}
\label{alphaheavydefinition}
\end{definition}
\begin{remark}\label{initialalphaheavywords}
 Trivially, if $A$ is $\A$-heavy (or light), then so is the initial factor $A_{0,j}$ for all $0 \leq j \leq m$.
\end{remark}
The following lemma is not difficult, but will be of great use in \S\ref{sturmian section}:
\begin{lemma}[The Reversing Principle]{Assume that $A$ is of length $n+1$, such that $A=a_0 \ldots a_{n-1}$ is $\A$-heavy, but \[\overline{w}(a_0 \ldots a_{n-1} a_n) \leq \A.\]  Then the word $(a_0 \ldots a_n)^{T}=a_n a_{n-1} \ldots a_0$ is $\A$-light.  Equivalently, the word \break $\rho(A_{0,n+1})$ is $(-\A)$-heavy.}\label{reversing lemma}
\begin{proof}
 Assume that there is some $i \in \{0,1,\ldots,n\}$ such that $\overline{w}(a_n a_{n-1} \ldots a_{n-i})>\A$.  We may assume $i \neq n$.  Then as $\overline{w}(a_0 \ldots a_{n-i}) \geq \A$ (the word $a_0 \ldots a_{n-1}$ was $\A$-heavy), we clearly have $\overline{w}(a_0 \ldots a_n)>\A$, a contradiction.
\end{proof}
\end{lemma}

>From this point, we refrain from statements in terms of both lightness and heaviness; we refer only to heaviness properties, but analogous statements regarding lightness are all possible.  Recall that a set $\mathcal{A} \subset \R$ is well-ordered by `$\geq$' if it contains no infinite increasing sequence.  We take `$\geq$' to be our standard ordering; for lightness, `$\leq$' would be the relevant ordering.

\begin{lemma}
 \label{inductive well ordered blocks}{Let $\mathcal{A} \subset \R$ be well-ordered.  Then for every $n \in \N$, the set $B_n=\{\overline{w}(A) \ST A \in \mathcal{A}^n\}$ is also well-ordered.}
 \begin{proof}
   It suffices to prove that if both $\mathcal{A}$ and $\mathcal{B}$ are well-ordered, then so is the set $(\mathcal{A}+\mathcal{B})=\{\gamma \ST \gamma= \alpha+\beta, \alpha \in \mathcal{A}, \beta \in \mathcal{B}\}$: our $B_n$ are subsets of $(\mathcal{A}+\mathcal{A}+\ldots+\mathcal{A})/n$, which is well-ordered if and only if $(\mathcal{A}+\mathcal{A}+\ldots+\mathcal{A})$ is well-ordered.  To the contrary, let $\gamma_1, \gamma_2, \ldots$ be an increasing sequence in $\mathcal{A}+\mathcal{B}$, where $\gamma_i=\alpha_i+\beta_i$, $\alpha_i \in \mathcal{A}$, $\beta_i \in \mathcal{B}$.  As $\mathcal{A}$ is well-ordered, we may pass to a subsequence $\gamma_{n(i)}$ such that $\alpha_{n(i)}$ is monotonically decreasing (not necessarily strictly): define \[n(i+1)= \min\{n >n(i) \ST \alpha_n=\max\{\alpha_j \ST j \geq n(i)\}\} \] and see that $\alpha_{n(i+1)} \leq \alpha_{n(i)}$, and each $n(i)$ is defined by well-orderedness of $\mathcal{A}$.  As $\gamma_{n(i)}$ are increasing, and the $\alpha_{n(i)}$ are nonincreasing, the $\beta_{n(i)}$ must be an infinite increasing sequence in $\mathcal{B}$, a contradiction.
 \end{proof}
\end{lemma}

We now investigate just how frequently one may expect to find $\alpha$-heavy factors of arbitrary sequences, depending on the target value $\alpha$.

\begin{lemma}{Let $\A > -\infty$ and $X$ be a sequence such that: \[\liminf_{n \rightarrow \infty} \overline{w}(X_{0,n})=\A.\]  Then for any $\delta<\A$, there is an $N \in \N$ such that the sequence $X_{N,\infty}=x_N x_{N+1}\ldots$ is $\delta$-heavy.}\label{boshrisingsun}
\begin{proof}
 Fix $\delta<\alpha$, and assume to the contrary that for every $N$, there is some $f(N)$ such that $\overline{w}(X_{N,f(N)})<\delta$.  Set $k_0=0$ and recursively define $k_i=f(k_{i-1})$.  Then represent $X=X_{k_0,f(k_0)}X_{k_1,f(k_1)}\ldots$.  It is seen that for all $i$, $\overline{w}(X_{1,f(k_i)})<\delta<\alpha$, contrary to our assumption that $\liminf_{n \rightarrow \infty} \overline{w}(X_{0,n})=\A$.
\end{proof}
\end{lemma}
\begin{corollary}{Let $\{\Omega,\mu,T\}$ be a probability measure preserving system which is not ergodic, and let $f \in L^1(\Omega,\mu)$ be such that $f_*(\omega)$ is not almost everywhere equal to a constant. Then $\mu\left(\mathcal{H}(\N)\right)>0$.}\label{nonergodic corollary}
\begin{proof}
As $\int_{\Omega}f_* d\mu=\int_{\Omega}f d\mu$, let $\Gamma=\{\omega \ST f_*(\omega)>\int_{\Omega}f d\mu\}$.  By assumption, $\mu(\Gamma)>0$, and $\forall \omega \in \Gamma$, $\exists N$ such that $T^N(\omega) \in \mathcal{H}(\N)$ (by Lemma \ref{boshrisingsun}).  As $\Gamma$ is covered by the preimages of $\mathcal{H}(\N)$, $\mu(\Gamma)>0$, and $T$ preserves $\mu$, we must have that $\mu\left(\mathcal{H}(\N)\right)>0$.
 \end{proof}
\end{corollary}

It is not difficult to construct a sequence $X$ of finite \emph{upper density}: \[\limsup_{n \rightarrow \infty} \overline{w}(X_{0,n})=\A < \infty,\] such that $X$ does not have arbitrarily long $\A$-heavy factors (for example, \break \mbox{$x_n/(n+1)$} will construct a sequence of upper density one with no $1$-heavy factors whatsoever), but the following lemma will extend the idea of Lemma \ref{boshrisingsun} as far as possible:

\begin{definition}\label{upper banach density}Define the \emph{upper Banach density of $X$}, for $X \in \mathcal{A}^{\N}$ or $\mathcal{A}^{\Z}$ by \[d_B^*(X)=\limsup_{(b_i-a_i) \rightarrow \infty} \overline{w}(X_{a_i,b_i}).\]\end{definition}
\begin{theorem}{ The alphabet $\mathcal{A}$ has the property that every $X \in \mathcal{A}^{\N}$ contains arbitrarily long $d_B^*(X)$-heavy words if and only if $\mathcal{A}$ is well-ordered. }\label{wellorderedheavy}
\begin{proof}
First, assume $\mathcal{A}$ is well-ordered.  Then for any $X \in \mathcal{A}^{\N}$, $d_B^*(X)<\infty$.  If $d_B^*(X)=-\infty$, there is nothing to prove, so assume $d_B^*(X) \in \R$, and fix some $\delta<d_B^*(X)$.  Assume that there is some $N<\infty$ so that for every $i$, there is some $f(i)<N$ for which $\overline{w}(X_{i,i+f(i)})<\delta$.  Similarly to Lemma \ref{boshrisingsun}, represent $X$ as a string of concatenated words of average weight strictly less than $\delta$, but note that there is now have a universal bound on the length of the words.  It follows that $d_B^*(X) \leq \delta$: for very large $b_i-a_i$, words of length $b_i-a_i$ may be considered as a concatenation of factors of length no larger than $N$, of average weight less than $\delta$, plus small extra pieces at the end of bounded length and weight.  So, $X$ must contain arbitrarily long $\delta$-heavy factors for arbitrary $\delta<d_B^*(X)$.

Now assume that there is an $X$ and a bound $N$ on the length of any $d_B^*(X)$-heavy factors of $X$.  For any $\epsilon>0$, let $s=s(\epsilon)$ be a factor of length $N$ which is $(d_B^*(X)-\epsilon)$-heavy (but by assumption, not $d_B^*(X)$-heavy).  Define a decreasing sequence $\epsilon_i$ by fixing an arbitrary $\epsilon_0>0$, and defining: \[\epsilon_{i+1}=\frac{d_B^*(X)-\overline{w}(S(\epsilon_i))}{2}>0\]

Continuing this process, create a sequence of words $\{S(\epsilon_i)\}_{i=0}^{\infty}$ of length $N$ whose average weights are strictly increasing.  By Lemma \ref{wellorderedheavy}, there is a contradiction.  Therefore, $X$ must contain arbitrarily long $d_B^*(X)$-heavy factors.

The proof of the converse is much shorter: let $\{\A_i\}_{i=0}^{\infty}$ be a sequence in $\mathcal{A}$ which is strictly increasing.  Then $X=\A_0 \A_1 \ldots$ does not contain any $d_B^*(X)$-heavy factors of any length.
\end{proof}
\end{theorem}
\begin{corollary}{Let $X \in \mathcal{A}^{\N}$, where $\mathcal{A}$ is a well-ordered and compact subset of $\R$.  Then under the transformation $\sigma$, the space $\overline{O^+(X)}$ contains some $X'$ which is $d_B^*(X)$-heavy.}\label{banachheavywords}
\begin{proof}
 Consider any limiting sequence $X'$ of the words $x(n)$, $d_B^*(X)$-heavy length $n$ factors of $X$ (these $x(n)$ exist by Theorem \ref{wellorderedheavy}).  By construction, $d_B^*(X')=d_B^*(X)$, and \[\inf_{n \in \N} \overline{w}(X'_{0,n})=\limsup_{n \rightarrow \infty} \overline{w}(X'_{0,n}).\qedhere\]
\end{proof}
\end{corollary}

\subsection{Local Heaviness}\label{local heaviness section}
In the following section, we introduce a version of heaviness which does not depend on an arbitrary constant $\A$:
\begin{definition}{A word $A=a_0 a_1 \ldots a_{n-1}$ is said to be \emph{heavy}, or \emph{locally heavy}, (\emph{light}, or \emph{locally light}), if for all $1 \leq i \leq n$: \[\overline{w}(A_{0,i}) \geq \overline{w}(A) \hspace{.2 in} (\overline{w}(A_{0,i})\leq \overline{w}(A)).\]  Equivalently, for all $1 \leq i \leq n-1$: \[\overline{w}(A_{0,i}) \geq \overline{w}(A_{i,n}) \hspace{.2 in} (\overline{w}(A_{0,i})\leq \overline{w}(A_{i,n})).\] A sequence $X$ is heavy (light) if for every $i \geq 0$: \[\overline{w}(X_{0,i})\geq \limsup_{n \rightarrow \infty} \overline{w}(X_{0,n}) \hspace{.2 in} \left( \overline{w}(X_{0,i}) \leq \liminf_{n \rightarrow \infty} \overline{w}(X_{0,n})\right).\]}
\label{localheavydefinition}
\end{definition}
The word $A$ is locally heavy if and only if $A$ is $\overline{w}(A)$-heavy.  However, the `target value' in this case \emph{varies with the word in question}, whereas in definition \ref{alphaheavydefinition} there was a preordained $\A$.  In the case $\mathcal{A}=\{0,1\}$, light words are called \emph{Lyndon words}, an object of study in combinatorics and computer science (see \cite{lothaire}).  Again, however, we will suppress statements regarding light words.

\begin{remark}
 In contrast to $\A$-heaviness, it is generally the case that for a given heavy $A$ of length $n$, there may be some $1 < i < n$ for which $A_{0,i}$ are not be locally heavy.  Consider the word $1010$, which is heavy, and the initial factor $101$, which is not.
\end{remark}

\begin{lemma}{Let $A$ and $B$ be heavy words.  Then the concatenation $AB$ is (locally) heavy if and only if $\overline{w}(A) \geq \overline{w}(B)$.}\label{concatenateheavy}
\begin{proof}
 The necessity is obvious: let $\abs{A}$ and $\abs{B}=m$.  If $\overline{w}(A)<\overline{w}(B)$, then: \[\overline{w}\left((AB)_{0,n}\right)=\overline{w}(A)<\overline{w}(B)=\overline{w}\left((AB)_{n,m+n}\right).\]
 \noindent Now, assuming $\overline{w}(A)\geq \overline{w}(B)$, we see that for $i\leq n$: \[\overline{w}\left((AB)_{0,i}\right)\geq \overline{w}(A)\geq\overline{w}(AB)\] and for $i>n$: \[ \overline{w}\left((AB)_{0,i}\right)\geq \overline{w}(B_{0,i-n})\geq \overline{w}(B_{i-n,m})= \overline{w}\left((AB)_{i,n+m}\right).\qedhere \]
\end{proof}
\end{lemma}

\begin{theorem}{Fix an alphabet $\mathcal{A}$.  Then every $X \in \mathcal{A}^{\N}$ contains arbitrarily long heavy factors if and only if $\mathcal{A}$ is well-ordered.}
\begin{proof}
 First, assume that $\mathcal{A}$ is well-ordered.  Let $X\in \mathcal{A}^{\N}$ and $N<\infty$ be such that $X$ contains no heavy factors of length longer than $N$.  Then represent $X$ as a chain of heavy words in the following manner:

 Let $A_1=x_0 \ldots x_{n_1-1}$ be the longest possible heavy factor beginning at $x_0$.  By assumption, $n_1\leq N$.  Let $A_2=x_{n_1} \ldots x_{n_1+n_2-1}$ be the longest heavy factor beginning at $x_{n_1}$ and again note that $n_2 \leq N$.  Continue in this manner to write $X=A_1 A_2 A_3 \ldots$ where each $\abs{A_i} \leq N$.  In light of Lemma \ref{concatenateheavy}, the average weights of these blocks must be strictly increasing.  Furthermore, because the length of each $A_n$ is bounded, there must be some specific length which occurs infinitely often, so there is an infinite collection of words of the same length, with strictly increasing average weight.  By Lemma \ref{inductive well ordered blocks}, this is impossible.

 Now, suppose that $\mathcal{A}$ has an infinite subsequence $\{\alpha_i\}_{i=0}^{\infty}$ which is strictly increasing.  Then the sequence $X=\alpha_0 \alpha_1 \ldots$ is seen to have no heavy words of length longer than one.
\end{proof}
\end{theorem}

\begin{remark}
 We do not claim that every $X$ has heavy factors of every length!  Consider the sequence $101010\ldots$; the alphabet $\{0,1\}$ is certainly well-ordered, but $X$ does not have any heavy factors of odd length larger than one.
\end{remark}

\section{The Morse-Thue sequence and substitution systems}\label{section-substitution}
In this section, we will define the classical Morse sequence, using it as an example to discuss certain aspects of heaviness.  After discussing this sequence, we will make brief remarks extending these properties to a general class of substitution systems.

\subsection{Heaviness in the Classical Morse(-Thue-Prouhet) Sequence}
\label{classicalmorse}

The classical Morse sequence may be built in the following manner.  Let $M(0)=0$, and $M(i+1)=M(i)\overline{M(i)}$.  The sequence $M$ such that $M_{0,2^n}=M(n)$ is the one-sided Morse sequence:
\[M(0)=0,\hspace{.02 in}M(1)=(0)(\overline{0})=01 , \hspace{.02 in}M(2)= (01)(\overline{01})=0110,\hspace{.02 in} M(3)=(0110)(\overline{0110}), \ldots\]
and the two-sided Morse sequence is given by the word $\hat{M} = M^T M$, centered about $\hat{M}(0)=M(0)$ (the decimal point appears to the left of $M_0$): \[M=.0110100110010110\ldots\] \[\hat{M}=\ldots 0110010110.0110100110\ldots.\]

For a survey of this history of this interesting sequence, including numerous applications and information on the many independent formulations, see \cite{allouche-shallit}.  It is easily seen that $M_{2i}M_{2i+1}\in \{01,10\}$ (including $i < 0$ in the case of $\hat{M}$).

\begin{lemma}{Let $A$ be a word of length $2k< n \leq 2k+2$ which is a factor of either $M$ or $\hat{M}$.  Then $k\leq w(A) \leq k+2$.}
 \begin{proof}
 For a word of the form $A'=M_{2i,2(i+k)}$, $w(A')=k$.  Then $A$ can only be of the form $A'$, $M_{2i-1}A'$, $A'M_{2(i+k)}$, or $M_{2i-1}A'M_{2(i+k)}$.  By considering all choices ($0$ and $1$) for values of $M_{2i-1}$ and $M_{2(i+k)}$, establish the inequality.
 \end{proof}
\end{lemma}

\begin{corollary}{Let $X$ be any sequence which is a factor of the Morse sequence such that $X_{0,2}=11$.  Then $X$ is $\frac{1}{2}$-heavy.  Similarly, if $X$ is a sequence which is a factor of the Morse sequence, and $X$ begins with $00$, then $X$ is $\frac{1}{2}$-light.}
\begin{proof} Pick an initial word $X_{0,i}$ of the form $11X'$ where $2k< \abs{X'}\leq 2k+2$.  Then $w(X_{0,i})=2+w(X')\geq k+2$, so $\overline{w}(X_{0,i})\geq\frac{1}{2}$.  The proof is similar for $X_{0,2}=00$.\end{proof}
\end{corollary}

We now appeal to a well-known result: the system generated by $M$ is uniquely ergodic (see \cite[Ch. 5]{luminysubstitutions}).  The unique invariant measure $\mu$ assigns $\mu(\Gamma)=\frac{1}{2}$, where $\Gamma=\{\omega \in \overline{O^+(M)} \ST \omega_0=1\}$, and $\mu(\Delta)>0$, where $\Delta=\{\omega \ST \omega_0 \omega_1=11\}$ (in fact, $\mu(\Delta)=\frac{1}{6}$).
\begin{corollary}{In the system $\{\overline{O^+(\hat{M})}, \mu, \sigma\}$, with the function $\chi_{\Gamma}$, where \break \mbox{$\Gamma=\{ \omega \ST \omega_0=1\}$,} we have $\mu \left( \mathcal{H}^{\chi_{\Gamma}}_{\sigma}(\N) \right) \geq \mu \left( \mathcal{H}^{\chi_{\Gamma}}_{\sigma}(\Z) \right) >0$, where $\hat{M}$ is the two-sided Morse sequence.  In the one-sided sequence space generated by $M$, $\mu \left(\mathcal{H}^{\chi_{\Gamma}}_{\sigma}(\N)\right)>0$.}
 \begin{proof}
    We have seen that $\Delta \subset \mathcal{H}_{\sigma}^{\chi_{\Gamma}}(\N)$, where $\Delta=\{\omega \ST \omega_0 \omega_1=11\}$, and we have already seen that $\mu (\Delta) >0$. Therefore, $\mathcal{H}_{\sigma}^{\chi_{\Gamma}}(\N)$ is of positive measure.  Similarly, in the two-sided sequence, $\mathcal{H}_{\sigma}^{\chi_{\Gamma}}(-\N)=\bigcap_{i=1}^{\infty}\mathcal{H}_{\sigma}^{\chi_{\Gamma}}(-i,0)$ is of positive measure; infinite words which end in $00$.  Therefore, in the two-sided sequence, any bisequence $X$ with $x_{-2}x_{-1}x_0 x_1=0011$ will be in $\mathcal{H}_{\sigma}^{\chi_{\Gamma}}(\Z)$, and the set $\Delta'=\{X \ST x_{-2}x_{-1}x_0 x_1=0011\}$ set is also seen to be of positive measure (to be precise, $\mu(\Delta ')=\frac{1}{12}$).
 \end{proof}
\end{corollary}
\subsection{General (Non-Mixing) Substitution Systems}\label{substitution systems}

We begin this subsection with the following theorem:
\begin{theorem}[Hal\'{a}sz \cite{halasz}]{Let $\{\Omega, T, \mu\}$ be a probability-measure-preserving system.  The value $E(\A)=e^{2 i \pi \A}$ belongs to the spectrum of $T$ if and only if there exists a set $\Gamma$ of measure $\A$ such that for $\mu$-almost every $\omega \in \Omega$, for all $n \in \N$: \begin{equation}\label{halasz equation}\abs{S_n(\omega)-n\A}\leq 1,\end{equation} where $f(\omega)=\chi_A(\omega)$.}\label{halasz}
\end{theorem}

\begin{corollary}{Let $\{\Omega, \mu, T\}$ be a probability-measure-preserving system, with $\A \in \Q$ such that $E(\A)$ is in the spectrum of $T$.  Then, letting $f(\omega)=\chi_{\Gamma}(\omega)$, where $\Gamma$ is a set which satisfies \eqref{halasz equation}, $\mu \left(\mathcal{H}(\N)\right)>0$.}\label{rational eigenvalue heavy sets}
\begin{proof}
 The quantities $S_n(\omega)-n\A$ are bounded almost everywhere (from Theorem \ref{halasz}) and discrete ($S_n(\omega) \in \Z$, and $\A \in \Q$).  Therefore, the minimum value of the sequence $\{S_n(\omega)-n\A\}_{n=0}^{\infty}$ is achieved at some minimal time $N(\omega)$ for almost all $\omega \in \Omega$.  So, $\omega \in T^{-N(\omega)}\left(\mathcal{H}(\N)\right)$.  As \[ \mu \left(\Omega \setminus \bigcup_{i=0}^{\infty}T^{-i}\left( \mathcal{H}(\N)\right) \right)=0,\] we must have $\mu \left( \mathcal{H}(\N)\right) >0$.
\end{proof}
\end{corollary}

We will provide a brief overview of a class of systems in symbolic dynamics with rational eigenvalues, as well as illustrating why satisfying \eqref{halasz equation} with an irrational eigenvalue does not guarantee positive measure heavy sets.

The Morse sequence may also be viewed as a fixed point of the substitution defined by $0 \rightarrow 01$ and $1 \rightarrow 10$:
 \[0 \rightarrow (01)=01 \rightarrow (01)(10)=0110 \rightarrow (01)(10)(10)(01)=01101001 \rightarrow \ldots \]

 In general, define a substitution system $\Sigma$ on $\Omega_k=\{0,1,\ldots,k-1\}$ by assigning $\Sigma(a) \in \Omega_k^{n(a)}$, for all $a \in \Omega_k$, and $n(a)<\infty$ for all $a$ (that is, $\Sigma$ assigns a word to each letter).  The substitution matrix for $\Sigma$ is: $A_{i,j}=\#\{j \textrm{ in } \Sigma(i)\}$.  In the case of the Morse sequence, the matrix is given by $A_{i,j} \equiv 1$.

Two different substitution systems might have the same matrix.  However, if the matrix $A$ is primitive ($\exists n \ST A^n_{i,j}>0 \hspace{.1 in}\forall i,j$), then the shift map defines a uniquely ergodic system $\{\Omega,\sigma\}$ on some limiting sequence $X \in \overline{\{\Sigma^n(x_0)\}_{n=1}^{\infty}}$ such that $\Sigma^N(X)=X$ for some $N \in \N$ (see \cite[Ch. 5]{luminysubstitutions}).

In the event that $\sigma$ is a substitution of \emph{constant length} ($n(a)$ is constant over $\Omega_k$), and $X$ is a periodic point under the substitution $\Sigma$, then the system $\{\overline{O^+(X)},\mu,\sigma\}$ has nontrivial rational eigenvalues (see \cite[Ch. 7]{luminysubstitutions}), and therefore there are nontrivial $\mu$-integrable functions $f:\overline{O^+(X)} \rightarrow \R$ with positive-measure heavy sets (in light of Corollary \ref{rational eigenvalue heavy sets}).

\begin{example}{Fix $\mathcal{A}=\{0,1,2\}$ and $\Sigma(0)=120$, $\Sigma(1)=201$, $\Sigma(2)=210$ (note that each $\Sigma(a)$ is a $1$-heavy word of length $3$ and average weight $1$).  Let $X$ be the invariant limiting sequence $\lim_{n \rightarrow \infty }\Sigma^n(2)=210201120210120\ldots$, and create a uniquely ergodic shift system $\{\Omega, \mu, \sigma\}$ (the substitution is primitive).  If we are interested in making the system invertible, we may create the natural invertible extension, as outlined in \cite[pp 239-241]{cornfeld-sinai}.  Define $f(\omega)=x_0$.  The following progression is extremely similar to that carried out in \S\ref{classicalmorse}, and details are omitted.

First, to establish a `target,' compute: $\int_{\Omega}f d\mu=1$.  Now, because $X_{3k,3(k+i)}$ are always of average weight $1$, consideration of blocks which may precede or follow these `evenly weighted' blocks (for example, $x_{3k+1}x_{3k+2}\neq 21$) yields that $\abs{S_n(X)-n} \leq 1$ for any sequence $X \in \Omega$.  So, any factor which begins with $2$, $12$, or $112$ will be heavy, and if $A$ is a factor ending with $0$, $01$, or $011$ then $\rho(A)$ is $-1$-heavy:
 \[\mathcal{H}(\N)\supset\{X \ST x_0 =2, \hspace{.1 in}x_0 x_1 =12, \hspace{.1 in} x_0 x_1 x_2 =112\},\]
 \[\mathcal{H}(\Z)\supset \{X \ST x_{-2}x_{-1}x_0 x_1 \in \{1012,1021,0112,0121\}\}.\]
Finally, quick density computations verify that $\mu\left( \mathcal{H}(\N)\right)\geq \mu \left(\mathcal{H}(\Z)\right)>0$.
}
\end{example}

To see the importance of rational eigenvalues to produce positive-measure heavy sets, let $\A \notin \Q$ and consider the circle rotation $\{S^1,\mu,R_{\A}\}$, where $\mu$ is Lebesgue measure and $R_{\A}(\omega)=\omega+\A \mod{1}$.  Fix $\Gamma=[0,\A)$ and $f(\omega)=\chi_{\Gamma}(\omega)$.  Then $\abs{S_n(\omega)-n\A}$ is bounded \cite{kesten1}.  However, as $\A \notin \Q$, these $S_n(\omega)-n\A$ do not necessarily ever achieve their infimum.  Indeed, we will see in Corollary \ref{uniquesturmian} that $\mathcal{H}_{R_{\A}}^f(\N)$ is exactly one point, and therefore a null set.

\section{Sturmian sequences}\label{section-sturmian}

\begin{definition}\label{sturmian definition}{A binary word, sequence, or bi-sequence $X$ is called \emph{Sturmian} if factors of the same length differ in weight by at most one.}
\end{definition}
\begin{remark}In a Sturmian sequence, the \emph{density} of the sequence exists \cite[Ch. 6]{luminysubstitutions}: \[d(X) = \lim_{n \rightarrow \infty} \overline{w}(X_{0,n}).\]  We restrict our attention to those Sturmian $X$ of density $\A \notin \Q$; Sturmian sequences of rational density are eventually periodic.
\end{remark}
\subsection{$\A$-Heaviness in Sturmian Sequences}\label{sturmian section}

We note the following theorem:
\begin{theorem}[E. Coven, G. Hedlund \cite{coven-hedlund}]{A bi-infinite binary sequence $X$ of minimal complexity ($p(n)+1$) is Sturmian if and only if for every $A$ which is a factor of $X$, $A^T$ is also a factor of $X$.} \label{coven-hedlund}
\end{theorem}
\noindent and use it to derive the following:
\begin{theorem}{A bi-infinite binary sequence $X$ is Sturmian of irrational density $\A$ if and only if for any $n \in \N$, there is a unique factor $A$ of length $n$ such that $A$ is $\alpha$-heavy, and a unique factor $A$ of length $n$ such that $\rho(A)$ is $(-\A)$-heavy.}\label{theorem-sturmian equivalence}
\begin{proof}
 Assume $X$ is Sturmian, and let $\A \notin \Q$ be its density.  Then by Theorem \ref{wellorderedheavy}, $X$ contains arbitrarily long $\A$-heavy factors ($\mathcal{A}=\{0,1\}$ is certainly well-ordered).  It is also seen that $X$ contains as factors arbitrarily long transposes of $\A$-light factors (using the `$\alpha$-light' version of Theorem \ref{wellorderedheavy}, and noting that $x_{-1}x_{-2}x_{-3}\ldots$ is also Sturmian of density $\A$).  So, $X$ has arbitrarily long factors $A$ such that $\rho(A)$ is $(-\A)$-heavy.  Therefore, there exists at least one word of each length which satisfies our criteria.

 Let $A$, $B$ be two distinct factors of length $n$ which are $\A$-heavy.  Assume they are of minimal length $n>0$, so that $A=C0$ and $B=C1$ for some factor $C$.  Then $w(A)$ and $w(B)$ are the two possible weights for factors of length $n$ in $X$, and both weights are at least as large as $n\A$.  As $\A \notin \Q$, both are strictly larger than $n\A$.  Therefore, every factor of length $n$ has average weight larger than some $\A+\epsilon$, contradicting the fact that $X$ was of density $\A$.  The proof is similar for two factors of length $n>0$ whose reversals are $(-\A)$-heavy.  If $\A \in \Q$, the result does not hold: consider $X=\ldots(10)(10)1(10)(10)(10)\ldots$.  This sequence is Sturmian, and the factors $10$ and $11$ are both $1/2$-heavy.

 For the converse, it will suffice, in light of Theorem \ref{coven-hedlund}, to show that $X$ is of minimal complexity and all transposes of factors are also factors.  Assume that the following conditions all hold for $1 \leq i \leq n-1$ (they are easy to verify for $i=1$):

 \begin{enumerate}
   \item There is a unique $\alpha$-heavy factor of length $i$, and a unique factor of length $i$ whose reversal is $(-\A)$-heavy.  \textit{This condition has been assumed for all $i$.}\label{uhb}\\
   \item If $A$ is a factor of length $i$, then $A^T$ is a factor. \label{t-prop}\\
   \item The factor $A$ of length $i-1$ is a right (left) special factor if and only if $1A^T$ is $\A$-heavy ($\rho(1A)$ is $-\A$-heavy). \label{ts-prop}\\
 \end{enumerate}

 By combining (\ref{uhb}) and (\ref{t-prop}), $X$ contains unique $\A$-light factors of length $i$, and unique words of length $i$ whose reversals are $(-\A)$-light.  Adding (\ref{ts-prop}), $p(n-1)$, as there is a unique right (or left) special factor for all $i\leq n-1$.  Establishing (\ref{t-prop}) and (\ref{ts-prop}) for all $n$, then, would ensure both minimal complexity and admissibility of transposes, sufficient to show that $X$ is Sturmian.

 Assume, then, that there is a factor $A$ of length $n \geq 2$ such that $A^T$ is not a factor.  Then, by our inductive hypothesis, $A_{0,n-1}^T$ and $A_{1,n}^T$ are both factors.  As $X$ is a bisequence, both of these words have precursors and successors (they do not begin or end $X$), so the words \[\overline{a_{n-1}}A_{0,n-1}^T=\overline{a_{n-1}}A_{1,n-1}^Ta_0 \hspace{.1 in} \textrm{and} \hspace{.1 in} A_{1,n}^T\overline{a_0}=a_{n-1} A_{1,n-1}^T \overline{a_0}\] are factors.  Note, then, that $A_{1,n-1}^T$ is both a left and right special factor.  Therefore, by (\ref{ts-prop}), $1A_{1,n-1}$ is the unique $\A$-heavy factor of length $n-1$, and $A_{1,n-1}1$ is the unique $\A$-light factor of length $1-n$.  Then \[\left(1A_{1,n-1}=(A_{1,n-1}1)^T\right) \Rightarrow \left(A_{1,n-1}^T=A_{1,n-1}\right).\]  As $A^T$ is not a factor, $A \neq A^T$, so $a_0 \neq a_{n-1}$.  Without loss of generality, let $a_0=1$ and $a_{n-1}=0$.  Let $B=B^T=A_{1,n-1}$, for convenience.

 The following have been shown to be factors of $X$: $1B0$, $0B0$, and $1B1$.  As $1B$ was $\A$-heavy, $1B1$ is the unique $\A$-heavy factor of length $n$.  Therefore, $1B0$ is not $\A$-heavy, and by the Reversing Principle (Lemma \ref{reversing lemma}), the word $\rho(1B0)$ is $(-\alpha)$-heavy.  Then certainly $\rho(0B0)$ is $(-\A)$-heavy and of the same length, contradicting (\ref{uhb}), our original assumption.  Therefore, (\ref{t-prop}) holds for factors of length $n$ as well.

 It remains only to show (\ref{ts-prop}) for factors of length $n$.  Let $A$ be the unique right special factor of length $n-1$; $A1$ and $A0$ both appear.  Inductively, then, $1A^T$ is the unique $\alpha$-heavy factor of length $n-1$.  However, by (\ref{t-prop}), both $1A^T$ and $0A^T$ appear.  So, if $1A^T$ is \emph{not} $\alpha$-heavy, it follows (again, by the Reversing Principle) that $A1$ and $A0$ are both $\alpha$-light, so that the reversals of the transposes are both $(-\A)$-heavy, contradicting (\ref{uhb}).

 In the reverse, let $1A$ be the unique $\alpha$-heavy factor of length $n$.  Then $1A_{0,n-2}$ is the unique $\A$-heavy factor of length $n-1$, so $A_{0,n-2}^T$ is a right special factor.  If $A^T$ is \emph{not}a right special factor, then $\overline{a_{n-2}}A_{1,n-2}^T$ is, and the previous reasoning would ensure that $1A_{1,n-2}\overline{a_{n-2}}$ would be $\alpha$-heavy, contradicting (\ref{uhb}).

 So, there is a $1-1$ correspondence between $\alpha$-heavy factors of length $n$ (which exist and are unique) and right special factors of length $n-1$.  Therefore, $p(n)+1$ for all $n$, and as $X$ contains as factors all transposes of factors, $X$ is Sturmian.
\end{proof}

\end{theorem}

\begin{remark}
 The heart of the above theorem is that Sturmian words can be characterized as having unique special factors, and when special factors are unique, they can be characterized by a heaviness condition.
\end{remark}

\begin{corollary}{Let $X$ be a Sturmian sequence (bisequence), of density $\A \notin \Q$.  Then $\overline{O^+(X)}$ contains exactly one sequence (bisequence) which is $\A$-heavy.}\label{uniquesturmian}
\begin{proof}
Let $x(n)$ be the $\A$-heavy factor of $X$ which is of length $n$.  It is seen that for $m>n$, $m,n \in \N$, $x(m)=x(n)A$ for some word $A$ of length $m-n$.  Therefore, in the compact space $\{0,1\}^{\N}$, let \[X'=\lim_{n \rightarrow \infty}x(n).\]  The sequence $X'$ is unique, is in the system generated by $X$, and $X'$ uniquely extends on the left as a Sturmian bisequence (see \cite[Ch. 6]{luminysubstitutions}).
\end{proof}
\end{corollary}
\subsection{Local Heaviness in Sturmian Sequences}
We will now approach the subject of Sturmian sequences using local heaviness (Definition \ref{localheavydefinition}), rather than $\A$-heaviness.
\begin{lemma}{There is exactly one heavy Sturmian word of length $n$ and weight $m$ for any choice $0 \leq m \leq n <\infty$.}\label{uniquesturmianheavy}
\begin{proof}
 For convenience, denote a word of length $n$ and weight $m$ as \emph{type} $(m,n)$, and the claim is apparent if $m=0$ or $n=1$.  Assume, then, that the claim is true for all words of length smaller than $n$.

 Let $A$ be Sturmian and heavy, of type $(pk_1, pk_2)$, where $k_1$ and $k_2$ are relatively prime.  Then applying the pigeonhole principle and Lemma \ref{concatenateheavy}, $A=A_1 A_2 \ldots A_p$, where each $A_i$ is Sturmian and heavy, of type $(k_1,k_2)$. Therefore, it is sufficient to prove the claim in the event when $m$ and $n$ are relatively prime.  If $m>n/2$, $A$ must be of the form $1^{n_0}01^{n_1}0\ldots 1^{n_k}0$, where the Sturmian condition requires that each $n_i=N$ or $N+1$ for some $N$.  As $(m,n)=1$, both values occur (if all $n_i=N$, for instance, $n=(k+1)(N+1)$ and $m=(k+1)N$), and the heaviness condition requires that $n_0=N+1$ and $n_k=N$.  Define $f(1^{n_i}0)_i-N$, and associate to $A$ the smaller word \[B=f(1^{n_0}0)f(1^{n_1}0)\ldots f(1^{n_k}0).\]

 We now show that $B$ is Sturmian.  Let $B'$ and $B''$ be factors of $B$ of equal length such that $w(B')=w(B'')+2$.  Then consider the two factors of $A$, $f(A')=B'$ and $f(A'')=B''$: $w(A')=w(A'')+2$, but $\abs{A'}=\abs{A''}+2$ as well.  However, by assuming the minimality on the length of $B'$ and $B''$, $B''$ begins with a zero.  Therefore, $A'' \neq A_{0,N}$; $A''$ is preceded by a zero.  Let $C_1 = 0 A''$.  Also, the last element of $A'$ must be a zero, so let $C_2=A'_{0,\abs{A'}-1}$.  Then $w(C_2)=w(C_1)+2$, but $\abs{C_1}=\abs{C_2}$, contradicting the assumption that $A$ is Sturmian.

 Now, to show that $B$ is heavy, begin with knowledge that $A$ is heavy and recall $n_i=N+f(1^{n_i}0)$:
 \begin{align*}
 \overline{w}(1^{n_0}0\ldots 1^{n_{i-1}}0) &\geq \overline{w}(1^{n_{i}}0\ldots 1^{n_{k-1}}0)\\
 \left(\sum_{j=0}^{i-1} n_j \right) \left( \sum_{j=i}^{k-1} (n_j+1)\right) &\geq \left(\sum_{j=i}^{k-1} n_j \right) \left( \sum_{j=0}^{i-1} (n_j+1) \right)\\
 \left(iN+w(B_{0,i}) \right) \left( (k-i)(N+1)+w(B_{i,k})\right) &\geq \left( (k-i)N+w(B_{i,k}\right) \left( i(N+1)+w(B_{0,i}) \right)\\
 (k-i)w(B_{0,i}) &\geq iw(B_{i,k})\\
 \overline{w}(B_{0,i}) &\geq \overline{w}(B_{i,k})
 \end{align*}
 So $B$ is a heavy Sturmian word of smaller length, and by the inductive hypothesis $B$ is unique, and therefore $A$ is unique.

 The proof works similarly for $m<n/2$, by considering the lengths of blocks of zeroes.\end{proof}
\end{lemma}

\begin{corollary}{Let $X$ be a Sturmian sequence which contains at least two ones and two zeroes.  Let $N=\max\{n \ST X \textrm{ contains }1^n\textrm{ or } 0^n \textrm{ as a factor}\}$, noting that our assumptions guarantee $N<\infty$ (and either $11$ or $00$ is not a factor $X$).  Then $X$ has exactly two distinct heavy factors of lengths $n\leq N$, and at most one heavy factor of all other lengths.}\label{2 heavy then 1 if sturmian}
\begin{proof}
 All factors of $X$ are Sturmian, and given a fixed length $n$, there are at most two weights possible for factors of length $n$.  It is therefore clear in light of Lemma \ref{uniquesturmianheavy} that for any $n$, there are at most two heavy factors of length $n$.  Assume that $d(X)\geq 1/2$, so that $N$ is the longest string of ones.  Then $1^N0$ (where $1^N$ represents a string of $N$ consecutive ones) is a factor, so $1^{N-i-1}0$ and $1^{N-i}$ for $i=0,\ldots,N-1$ are also factors, giving two heavy factors of lengths $1,\ldots,N$.  The proof is similar if $N$ is the length of the longest string of ones.

 If $X$ has two distinct heavy factors of length $n>N$, then as they are Sturmian heavy words and distinct, they must be of two different weights, and the weights therefore differ by one.  Let $B$ and $C$ be the two factors, and let $w(B)=w(C)+1$. That $n>N$ ensures that each of them begins with a one and ends with a zero; a heavy binary word containing both $1$ and $0$ must begin with a one and end with a zero.  Then it is seen that $w(B_{0,n-1})=w(B)=w(C)+1=w(C_{1,n})+2$, which contradicts that $X$ is Sturmian.
\end{proof}
\end{corollary}

\begin{corollary}
 {Let $X$ be a Sturmian sequence whose density is $\alpha \notin \Q$.  Then there is a unique locally heavy sequence in $\overline{O^+(X)}$.}
 \begin{proof}
   Recall that for a sequence to locally heavy is, by definition, the same as the word being $\limsup_{n \rightarrow \infty}\overline{w}(X_{0,n})$-heavy.  We appeal to Corollary \ref{uniquesturmian}.
 \end{proof}
\end{corollary}

We conclude with a remark on the construction of heavy and $\A$-heavy Sturmian words.  Given an $\A \in [0,1]$ and $n \in \N$, define $a_i=\integer{i\A}-\integer{(i-1)\A}$ for $i=0,1,\ldots,n-1$, where $\integer{\omega}=\max\{n \in \Z \ST n \leq \omega\}$.  Then $A=a_0 \ldots a_{n-1}$ is the unique Sturmian word of type $(w(A),\abs{A})$, and if $\A \notin \Q$, $A$ is the unique $\A$-heavy factor of the infinite Sturmian sequence of density $\A$.  The sequence $a_i$ is related to the \emph{spectrum} of $\A$.  See \cite{graham-lin-lin} and \cite{boshernitzan-fraenkel}.

\section*{Acknowledgements}
Many thanks to Michael Boshernitzan for ideas, helpful suggestions, and interesting problems.  Also, David Damanik, Michael Keane, and Karl Petersen have all offered advice, provided valuable references, and aided development of ideas herein.  The author also wishes to thank the referees for many comments on improving the readability of this paper and implementing standard terminology.

\medskip

Received February 15, 2008; revised August 28, 2008.

\medskip

\end{document}